\documentclass[10pt,openbib]{article}
\usepackage{amsmath}
\usepackage{amsfonts}
\usepackage{amssymb}
\usepackage{graphicx}
\newcommand{\R}{\mathbb{R}}

\newcommand{\dee}{\mathop{\! \,\mathrm{d} \!}\nolimits}

\newcommand{\onehalf}{\mbox{$\frac{\scriptstyle 1}{\scriptstyle 2}$}}

\newcommand{\ttfrac}[2]{\mbox{$\frac{{\scriptstyle #1}}{{\scriptstyle #2}}$}}

\allowdisplaybreaks

\begin{document}
\thispagestyle{empty}
\begin{center}
{\Large \textbf{Normalization and reduction of the \\
\rule{0pt}{16pt}Stark Hamiltonian}}
\mbox{}\vspace{.05in} \\ 
\mbox{}\\
Richard Cushman\footnotemark 
\end{center}
\footnotetext{e-mail: r.h.cushman@gmail.com  \hfill version: \today } \medskip

\begin{abstract}
This paper details a calculation of the second order normal form of the Stark effect Hamiltonian 
after regularization, using the Kustaanheimo-Stiefel mapping. After reduction we 
obtain an integrable two degree of freedom system on $S^2_h \times S^2_h$, which we reduce 
again to obtain a one degree of freedom Hamiltonian system. 
\end{abstract} \bigskip

Most papers on the Stark Hamiltonian discuss only the two degree of freedom version, 
which was shown to be completely integrable by Lagrange \cite{lagrange}, see also \cite{lantoine-russell}. 
We do not expect the three degree of freedom Stark Hamiltonian to be completely integrable. Our 
treatment follows that of \cite{cushman-sadovskii}, which does not give any details how the 
second order normal form was obtained. We use the notation and results of \cite{cushman22A}. 

\section{The basic setup}

On $T_0{\R }^3 = ({\R }^3\setminus \{ 0 \}) \times {\R }^3$ with coordinates $(x,y)$ and standard symplectic form ${\omega }_3= \sum^3_{i=1}\dee x_i \wedge \dee y_i$ consider the Stark Hamiltonian 
\begin{equation}
K(x,y) = \onehalf \langle y, y \rangle -\frac{1}{|x|} + fx_3.
\label{eq-one}
\end{equation}
Here $\langle \, \, , \, \, \rangle $ is the Euclidean inner product on ${\R }^3$ with associated norm 
$|\, \, |$. On the negative energy level $-\onehalf k^2$ with $k>0$ rescaling time $\dee t \mapsto 
\frac{|x|}{k}\dee s$, we obtain
\begin{equation}
0 = \frac{1}{2k}\big( |x| \langle y,y \rangle +k^2|x|\big) -\frac{1}{k} + fx_3 \frac{|x|}{k}. 
\label{eq-two}
\end{equation}
In other words, $(x,y)$ lies in the $\frac{1}{k}$ level set of 
\begin{equation}
\widehat{K}(x,y) = \frac{1}{2k}|x|\big( \langle y, y \rangle +k^2 |x| \big) + f x_3\frac{|x|}{k}. 
\label{eq-three}
\end{equation}
We assume that $f$ is small, namely, $f = \varepsilon \beta $. After the symplectic coordinate change $(x,y) \mapsto (\frac{1}{k}x, ky)$ the Hamiltonian $\widehat{K}$ becomes the preregularized Hamiltonian 
\begin{equation}
\mathcal{K}(x,y) = \onehalf |x|(\langle y, y \rangle +|x|) +\varepsilon \beta x_3|x|
\label{eq-four}
\end{equation}
on the level set ${\mathcal{K}}^{-1}(1)$. \medskip 

Let $T_0{\R }^4 = ({\R}^4\setminus \{ 0 \} ) \times {\R }^4$ have 
coordinates $(q,p)$ and a symplectic form ${\omega }_4 = \sum^4_{i=1}\dee q_i \wedge \dee p_i$. 
Pull back $\mathcal{K}$ by the Kustaanheimo-Stiefel 
mapping 
\begin{displaymath}
\mathcal{KS}: T_0{\R }^4 \rightarrow T_0{\R }^3: (q,p) \mapsto (x,y),
\end{displaymath}
where 
\begin{align*}
x_1 & = 2(q_1q_3+q_2q_4) = U_2-K_1 \\
x_2 & = 2(q_1q_4-q_2q_3) = U_3 -K_2 \\
x_3 & = q^2_1+q^2_2-q^2_3-q^2_4 = U_4 - K_3 \\
y_1 & = (\langle q,q\rangle )^{-1}(q_1p_3+q_2p_4 +q_3p_1 +q_4p_2) = (H_2 + V_1)^{-1}V_2 \\
y_2 & = (\langle q,q\rangle )^{-1}(q_1p_4-q_2p_3 -q_3p_2 +q_4p_1) = (H_2 + V_1)^{-1}V_3 \\
y_3 & = (\langle q,q\rangle )^{-1}(q_1p_1+q_2p_2 +q_3p_3 +q_4p_4) = (H_2 + V_1)^{-1}V_4 
\end{align*}
\noindent and 
\begin{align*}
H_2 & = \onehalf (p^2_1+p^2_2+p^2_3+p^2_4 +q^2_1+q^2_2 +q^2_3+q^2_4) \\
\Xi & = q_1p_2-q_2p_1 +q_3p_4 -q_4p_3,
\end{align*} 
to get the regularized Stark Hamiltonian 
\begin{equation}
\mathcal{H} = H_2 +\varepsilon \beta \big( U_4V_1+H_2U_4 -K_3V_1 -H_2K_3 \big)  
\label{eq-four}
\end{equation}
on ${\Xi}^{-1}(0)$, since $|x| = \langle q, q \rangle = H_1+V_1$. Here 
\begin{align}
K_1 & =-(q_1q_3+q_2q_4+p_1p_3 +p_2p_4) \notag \\
K_2 &= -(q_1q_4-q_2q_3+p_1p_4-p_2p_3)  \notag \\
K_3 & = \onehalf (q^2_3+q^2_4 +p^2_3+p^2_4 -q^2_1 -q^2_2 -p^2_1-p^2_2) \notag \\
\rule{0pt}{12pt} L_1 & = q_4p_1-q_3p_2+q_2p_3-q_1p_4 \notag \\
L_2 & = q_1p_3+q_2p_4-q_3p_1-q_4p_2 \notag  \\
L_3 & =  q_3p_4-q_4p_3+q_2p_1-q_1p_2 \notag \\
\rule{0pt}{12pt} U_1 & = -(q_1p_1+q_2p_2+ q_3p_3+ q_4p_4) \notag \\
U_2 & = q_1q_3+q_2q_4-p_1p_3-p_2p_4 \notag \\
U_3 & = q_1q_4 - q_2q_3 +p_2p_3 - p_1p_4 \notag \\
U_4 &  = \onehalf (q^2_1+q^2_2 -q^2_3 -q^2_4 +p^2_3 +p^2_4-p^2_1-p^2_2) \notag \\  
\rule{0pt}{12pt} V_1 & = \onehalf (q^2_1+q^2_2 +q^2_3+q^2_4-p^2_1 -p^2_2-p^2_3-p^2_4) \notag \\
V_2 & =  q_1p_3 +q_2p_4 +q_3p_1 + q_4p_2 \notag \\
V_3 & =  q_1p_4 - q_2p_3 +q_4p_1 - q_3p_2 \notag \\
V_4 & =  q_1p_1 + q_2p_2 -q_3p_3 -q_4p_4. \notag  
\end{align}
generate the algebra of polynomials invariant under the $S^1$ action 
${\varphi }^{\Xi}_s$ given by the flow of $X_{\Xi}$ on $(T{\R }^4 = {\R }^8, {\omega }_4)$. 
The Hamiltonian $\mathcal{H}$ (\ref{eq-four}) is invariant under this $S^1$ action and thus 
is a smooth function on the orbit space ${\Xi}^{-1}(0)/S^1 \subseteq {\R }^{16}$ with coordinates 
$(K,L,H, {\Xi};U,V)$.  

\section{The first order normal form on ${\Xi}^{-1}(0)/S^1$}

The harmonic oscillator vector field $X_{H_2}$ on $(T{\R }^4, {\omega }_4)$ induces the vector field 
$Y_{H_2} = \sum^4_{i=1}\big( 2V_i \frac{\partial }{\partial U_i} - 2U_i\frac{\partial }{\partial V_i }\big) $ on 
the orbit space ${\R }^8/S^1 \subseteq {\R }^{16}$, which leaves ${\Xi}^{-1}(0)/S^1$ invariant. \medskip 

We now compute the first order normal form of the Hamiltonian $\mathcal{H}$ (\ref{eq-four}) on the reduced space ${\Xi }^{-1}(0)/S^1 \subseteq {\R }^8/S^1$. \medskip 

The average of $H_2U_4-K_3V_1$ over the flow
\begin{displaymath}
{\varphi }^{Y_{H_2}}_t(K,L,H_2, \Xi ; U,V) = 
\big( K,L, H_2, \Xi ; U \cos 2t + V\sin 2t, -U\sin 2t + V \cos 2t \big) 
\end{displaymath}
of $Y_{H_2}$ is 
\begin{align*}
\overline{H_2U_4-K_3V_1} & = \ttfrac{1}{\pi }\int^{\pi }_0 ({\varphi }^{Y_{H_2}}_t)^{\ast }(H_2U_4 - K_3V_1 )\, 
\dee t \\
& \hspace{-.75in}= \ttfrac{1}{\pi }\int^{\pi }_0 (U_4\cos 2t +V_4 \sin 2t) \dee t \, H_2 -
        \ttfrac{1}{\pi }\int^{\pi }_0 (-U_1\sin 2t +V_1 \cos 2t) \dee t \, K_3 =0.  
\end{align*}
The second equality above follows because $L_{X_{H_2}}K_3 =0$ and the third because $\overline{\cos 2t} = \overline{\sin 2t } =0$. The average of $U_4V_1$ over the flow 
of $Y_{H_2}$ on ${\Xi }^{-1}(0)/S^1$ is 
\begin{align*}
\overline{U_4V_1} & = \ttfrac{1}{\pi }\int^{\pi }_0 ({\varphi }^{Y_{H_2}}_t)^{\ast }(U_4V_1 )\, \dee t \\
& = -\onehalf U_1U_4 \, \overline{\sin 4t} + U_4V_1\, \overline{{\cos }^2 \, 2t} - 
U_1V_4\, \overline{{\sin }^2 \, 2t} +\onehalf V_1V_4 \, \overline{\sin 4t} \\ 
& = \onehalf (U_4V_1 - U_1V_4) = -\onehalf H_2K_3, 
\end{align*}
since $\overline{{\cos }^2 \, 2t} = \overline{{\sin }^2 \, 2t} = \onehalf $ and 
$\overline{\sin 4t} =0$. The last equality above follows from the explicit 
description of the orbit space ${\R}^8/S^1$ as the 
semialgebraic variety in ${\R }^{16}$ with coodinates $(K,L,U,V;H_2, \Xi )$ given by
\begin{align}
\langle U,U \rangle = U^2_1+U^2_2+U^2_3 +U^2_4 & = H^2_2 - {\Xi }^2 \ge 0 \, \, \, H_2 \ge 0 \notag \\
\langle V,V \rangle = V^2_1+V^2_2+V^2_3 +V^2_4 & = H^2_2 - {\Xi }^2 \ge 0 \notag \\
\langle U,V \rangle = U_1V_1 +U_2V_2+U_3V_3 + U_4V_4 & =0 \notag \\
U_2V_1-U_1V_2 & = L_1\Xi - K_1H_2 \notag \\
U_3V_1-U_1V_3 & = L_2 \Xi - K_2H_2 
\label{eq-fivenwstar} \\
U_4V_1 - U_1V_4 & = L_3 \Xi  -K_3H_2 \notag \\
U_4V_3 -U_3V_4 & = K_1 \Xi -L_1H_2 \notag \\
U_2V_4 -U_4V_2 & = K_2 \Xi -L_2H_2 \notag \\
U_3V_2 -U_2V_3 & = K_3 \Xi -L_3H_2 . \notag 
\end{align}
So the average of $U_4V_1+ H_2U_4 - K_3V_1 - H_2K_3$ over the flow of $Y_{H_2}$ is 
$-\ttfrac{3}{2}H_2K_3$ on ${\Xi }^{-1}(0)/S^1$. Thus the first order normal form of the 
regularized Stark Hamiltonian $\mathcal{H}$ (\ref{eq-four}) on ${\Xi }^{-1}(0)/S^1$ is 
\begin{equation}
{\mathcal{H}}^{(1)}_{\mathrm{nf}} = H_2 - \ttfrac{3}{2} \beta \varepsilon H_2K_3. 
\label{eq-five}
\end{equation}

\section{The second order normal form on ${\Xi }^{-1}(0)/S^1$}

In order to compute the second order normal form of the Hamiltonian $\mathcal{H}$ on 
${\Xi }^{-1}(0)/S^1$, we need to find a function $F$ on ${\R }^{16}$ such that changing coordinates by the time $\varepsilon $ value of the flow of the Hamiltonian vector field $Y_F$ brings the regularized Hamiltonian 
$\mathcal{H}$ (\ref{eq-four}) into first order normal form. Choose $F$ so that  
\begin{equation}
L_{Y_F}H_2 = \beta \big( -U_4V_1 - \ttfrac{1}{2}H_2K_3 -H_2U_4 +K_3V_1 \big) .
\label{eq-seven}
\end{equation}
The following calculation shows that this choice does the job. 
\begin{align}
({\varphi }^{Y_F}_{\varepsilon })^{\ast }\mathcal{H} & = 
\mathcal{H} + \varepsilon L_{Y_F}\mathcal{H} + \onehalf {\varepsilon}^2 L^2_{Y_F}\mathcal{H} 
+ \mathrm{O}({\varepsilon}^3) \notag \\
& = H_2 +\varepsilon \beta \big( U_4V_1+H_2U_4-K_3V_1 -H_2K_3 \big) + \varepsilon L_{Y_F}H_2 
\notag \\
& \hspace{.25in} +{\varepsilon }^2\beta L_{Y_F}\big( U_4V_1+ H_2U_4 -K_3V_1 -H_2K_3 \big) 
+ \onehalf {\varepsilon }^2 L^2_{Y_F}H_2 + \mathrm{O}({\varepsilon }^3) 
\notag \\
& = H_2 +\varepsilon \beta \big( U_4V_1+H_2U_4-K_3V_1 -H_2K_3 \big) 
\notag \\
& \hspace{.25in} + \varepsilon \beta \big( -U_4V_1 - \ttfrac{1}{2}H_2K_3 -H_2U_4 +K_3V_1 \big) 
\notag \\
& \hspace{.5in} + {\varepsilon }^2 \big[ L_{Y_F}\big( -L_{Y_F}H_2 -\ttfrac{3}{2}\beta H_2K_3 \big) 
+\onehalf L^2_{Y_F}H_2 \big] + \mathrm{O}({\varepsilon }^3) 
\notag \\
& = H_2 -\ttfrac{3}{2}\varepsilon \beta H_2K_3 -\onehalf {\varepsilon }^2\big( L^2_{Y_F}H_2 
+3\beta L_{Y_F}(H_2K_3) \big) + \mathrm{O}({\varepsilon}^3). 
\label{eq-sevennwstar}
\end{align}

To determine the function $F$, we solve equation (\ref{eq-seven}). Write $F = F_1+F_2$, where  
$L_{Y_{H_2}}F_1 = \beta (U_4V_1+ \onehalf H_2K_3)$ and $L_{Y_{H_2}}F_2 = \beta (H_2U_4-K_3V_1)$. 
Then 
\begin{align}
L_{Y_F}H_2 & = -L_{Y_{H_2}}Y_F = -L_{Y_{H_2}}F_1 - L_{Y_{H_2}}F_2 \notag \\
& = -\beta (U_4V_1+ \onehalf H_2K_3) - \beta (H_2U_4-K_3V_1).
\label{eq-ninestar}
\end{align}
Since $L_{Y_{H_2}}V_4 = -2 U_4$ and $L_{Y_{H_2}}U_1 = 2V_1$, it follows that 
\begin{subequations}
\begin{equation}
F_2 = - \ttfrac{\beta }{2}(H_2V_4 +K_3U_1). 
\label{eq-eighta}
\end{equation}
Now 
\begin{align*}
F_1 & = \ttfrac{\beta }{\pi } \int^{\pi}_0 t({\varphi }^{Y_{H_2}}_t)^{\ast }\big( U_4V_1 + \onehalf H_2K_3 \big) 
\, \dee t = \ttfrac{\beta }{\pi } \int^{\pi}_0 t({\varphi }^{Y_{H_2}}_t)^{\ast }(U_4V_1) \dee t 
+ \ttfrac{\pi \beta }{4}H_2K_3, 
\end{align*}
see \cite{cushman}, and 
\begin{align*}
\ttfrac{\beta }{\pi } \int^{\pi}_0 t({\varphi }^{Y_{H_2}}_t)^{\ast }(U_4V_1) \dee t & = \\
&\hspace{-1.25in}=  -\ttfrac{\beta }{2} (U_1U_4) \, \ttfrac{1}{\pi } \int^{\pi }_0 t \, \sin 4t \dee t +
\beta (U_4V_1)\, \ttfrac{1}{\pi } \int^{\pi }_0 t \, {\cos }^2\, 2t \dee t \\
&\hspace{-.75in} -\beta (U_1V_4)\, \ttfrac{1}{\pi } \int^{\pi }_0 t \, {\sin }^2\, 2t \dee t 
+ \ttfrac{\beta }{2} (V_1V_4)\,  \ttfrac{1}{\pi } \int^{\pi }_0 t \, \sin 4t \dee t \\
&\hspace{-1.25in} = \ttfrac{\beta }{8} (U_1U_4 - V_1V_4) +\ttfrac{\pi \beta }{4}(U_4V_1 -U_1V_4) , 
\end{align*}
since $\frac{1}{\pi } \int^{\pi}_0 t \, \sin 4t \dee t = -\frac{1}{4}$ and 
$\frac{1}{\pi } \int^{\pi}_0 t \, {\sin }^2\, 2 t \dee t = \frac{1}{\pi } \int^{\pi}_0 t \,  {\cos }^2\, 2 t \dee t = 
\frac{\pi }{4}$. 
Thus 
\begin{equation}
F_1  = \frac{\beta }{8}(U_1U_4-V_1V_4) +\frac{\pi \beta }{4}(U_4V_1 -U_1V_4 + H_2K_3) \notag \\ 
=  \frac{\beta }{8}(U_1U_4-V_1V_4)   
\label{eq-eightb}
\end{equation}
on ${\Xi}^{-1}(0)/S^1$, see (\ref{eq-fivenwstar}).
\end{subequations}
Hence on ${\Xi}^{-1}(0)/S^1$ 
\begin{equation}
F = F_1+F_2 = \frac{\beta }{8}(U_1U_4-V_1V_4) -\frac{\beta }{2}(H_2V_4 + K_3U_1). 
\label{eq-nine}
\end{equation}

We now calculate the average over the flow of $Y_{H_2}$ of  
\begin{equation}
-\ttfrac{3}{2}\beta L_{Y_F}(H_2K_3) -\onehalf L^2_{Y_F}H_2,  
\label{eq-ten}
\end{equation}
which is the ${\varepsilon }^2$ term in the transformed Hamiltonian 
$({\varphi }^{Y_F}_{\varepsilon})^{\ast }\mathcal{H}$, see (\ref{eq-sevennwstar}). This determines 
the second order normal form of $\mathcal{H}$ on ${\Xi }^{-1}(0)/S^1$. We begin 
with the term 
\begin{displaymath}
-\ttfrac{3}{2} \beta L_{Y_F}(H_2K_3) = -\ttfrac{3}{2}\beta [K_3(L_{Y_F}H_2) -H_2(L_{Y_{K_3}}F)] .
\end{displaymath}
The average of 
\begin{displaymath}
-\ttfrac{3}{2}\beta \, K_3(L_{Y_F}H_2) = \ttfrac{3}{2}\beta \, K_3(U_4V_1 
+\onehalf \beta \, H_2K_3 + H_2U_4 -K_3V_1) 
\end{displaymath}
vanishes on ${\Xi}^{-1}(0)/S^1$. The term  
\begin{align*}
\ttfrac{3}{2}\beta H_2(L_{Y_{K_3}}F) & = \ttfrac{3}{2}{\beta }^2H_2\, L_{Y_{K_3}}\big( 
\ttfrac{1}{8}(U_1U_4-V_1V_4) -\onehalf (H_2V_4 +K_3U_1) \big) \\
&\hspace{-.75in}
 = \ttfrac{3}{2}{\beta }^2 H_2 \big[ \big( -2L_2\frac{\partial }{\partial K_1} +2L_1\frac{\partial }{\partial K_2} 
- 2K_2\frac{\partial }{\partial L_1} +2K_1\frac{\partial }{\partial L_2}  \\
&\hspace{-.5in} -2U_4\frac{\partial }{\partial U_1} + 2U_1 \frac{\partial }{\partial U_4} 
-2V_4\frac{\partial }{\partial V_1} + 2V_1\frac{\partial }{\partial V_4 } \big) \big]
\big( \ttfrac{1}{8}(U_1U_4-V_1V_4) \\
&\hspace{-.25in} -\onehalf (H_2V_4 +K_3U_1) \big) , \, \, \mbox{see \cite[table 1]{cushman22A} } \\
& \hspace{-.75in} = \ttfrac{3}{2}{\beta }^2H_2 \big[ \ttfrac{1}{4}(-U^2_4+U^2_1+V^2_4-V^2_1) 
-H_2V_1 +K_3U_4 \big] . 
\end{align*}
Next we calculate $\ttfrac{3}{2}\beta \overline{H_2(L_{Y_{K_3}}F)}$. Since $\overline{H_2V_1} = 0 
= \overline{K_3U_4}$ we need only calculate the average of $U^2_1$, $U^2_4$, $V^2_1$ and 
$V^2_4$. We get 
$\overline{U^2_1} = \onehalf (U^2_1+V^2_1) = \overline{V^2_1}$ and 
$\overline{U^2_4}  = \onehalf (U^2_4+V^2_4) = \overline{V^2_4} $. 
Thus $\ttfrac{3}{2}\beta \overline{H_2(L_{Y_{K_3}}F)} = 0$. 
So the average $-\frac{3}{2}\beta \overline{L_{Y_F}(H_2K_3)}$ of the first term in expression (\ref{eq-ten}) vanishes on ${\Xi}^{-1}(0)/S^1$. \medskip 

\vspace{-.15in}Next we calculate the average of the term $L^2_{Y_F}H_2$ in expression (\ref{eq-ten}) on 
${\Xi }^{-1}(0)/S^1$. We have
\begin{align*}
L^2_{Y_F}H_2 & = -L_{Y_F}(L_{Y_{H_2}}F)  \\
& = -\beta L_{Y_F}\big( U_4V_1 + \onehalf H_2K_3 +H_2U_4 -K_3V_1) , \, \, 
\mbox{using (\ref{eq-ninestar})} \\
& = \beta \big[ \overbrace{(L_{Y_{U_4}}F)V_1}^{I} + \overbrace{U_4(L_{Y_{V_1}}F)}^{II} 
-\onehalf \overbrace{(L_{Y_F}H_2)K_3}^{III} + \onehalf \overbrace{H_2(L_{Y_{K_3}}F)}^{IV} \\
& \hspace{.5in} -\overbrace{(L_{Y_F}H_2)U_4}^{V} + \overbrace{H_2(L_{U_4}F)}^{VI} - \overbrace{(L_{Y_{K_3}}F)V_1}^{VII} -\overbrace{K_3(L_{Y_{V_1}}F)}^{VIII} \big] . 
\end{align*}
We begin by finding 
\begin{align*}
L_{Y_{H_2}}F & = \beta \big[ \big( 2V_1 \frac{\partial }{\partial U_1} + 2V_2 \frac{\partial }{\partial U_2} 
+2V_3 \frac{\partial }{\partial U_3} + 2V_4 \frac{\partial }{\partial U_4} 
-2U_1 \frac{\partial }{\partial V_1} - 2U_2 \frac{\partial }{\partial V_2} \\
& \hspace{.25in} - 2U_3 \frac{\partial }{\partial V_3} -2U_4 \frac{\partial }{\partial V_4} \big) \big] \big( 
\ttfrac{1}{8}(U_1U_4-V_1V_4)  -\onehalf (H_2V_4 + K_3U_1) \big) \\ 
& = \beta \big[ \ttfrac{1}{2}(V_1U_4+U_1V_4) +H_2U_4 - K_3V_1 \big] ; \\
L_{Y_{K_3}}F & = \beta \big( \big( -2L_2 \frac{\partial }{\partial K_1} + 2L_1 \frac{\partial }{\partial K_2} 
-2K_2 \frac{\partial }{\partial L_1} + 2K_1 \frac{\partial }{\partial L_2} -2U_4 \frac{\partial }{\partial U_1} + 2U_1 \frac{\partial }{\partial U_4} \\
& \hspace{.25in} - 2V_4 \frac{\partial }{\partial V_1} +2V_1 \frac{\partial }{\partial V_4} \big) 
\big( \ttfrac{1}{8}(U_1U_4-V_1V_4) -\onehalf (H_2V_4+K_3U_1) \big) \big) \\
& = \beta \big[ \ttfrac{1}{4}(-U^2_4 +U^2_1 +V^2_4 - V^2_1) -H_2V_1 +K_3U_4 \big] ; \\
L_{Y_{U_4}}F & = \beta \big[ \big( -2U_1 \frac{\partial }{\partial K_1} - 2U_3 \frac{\partial }{\partial L_1} 
+2U_3 \frac{\partial }{\partial L_2} - 2V_4 \frac{\partial }{\partial H_2} 
-2K_3 \frac{\partial }{\partial U_1} + 2L_2 \frac{\partial }{\partial U_2} \\
& \hspace{.25in} + 2V_3 \frac{\partial }{\partial U_3} - 2H_2 \frac{\partial }{\partial V_4} \big) 
\big( \ttfrac{1}{8}(U_1U_4-V_1V_4) -\onehalf (H_2V_4+K_3U_1) \big) \big] \\
& = \beta \big[ V^2_4 +K^2_3 -\ttfrac{1}{4}(K_3U_4 - H_2V_1) +H^2_2 \big] ; \\
L_{Y_{V_1}}F & = \beta \big[ \big( -2V_2 \frac{\partial }{\partial K_1} - 2V_1 \frac{\partial }{\partial L_2} 
+2V_2 \frac{\partial }{\partial L_2} + 2U_2 \frac{\partial }{\partial H_2} 
+2H_2 \frac{\partial }{\partial U_2} + 2K_1 \frac{\partial }{\partial V_2} \\
& \hspace{.25in} + 2K_2 \frac{\partial }{\partial V_3} + 2U_4 \frac{\partial }{\partial V_4} \big) \big( \ttfrac{1}{8}(U_1U_4-V_1V_4) -\onehalf (H_2V_4+K_3U_1) \big) \big] \\
& = \beta (-U_2V_2 -\ttfrac{1}{4}U_4V_1) .
\end{align*}
\begin{subequations}
So the average of term I on ${\Xi}^{-1}(0)/S^1$ is
\begin{align}
\beta \overline{(L_{Y_{U_4}}F)V_1} & = {\beta }^2(\overline{V_1V^2_4} + \overline{K^2_3V_1} 
-\ttfrac{1}{4}\overline{K_3U_1V_4} + \ttfrac{1}{4}\overline{H_2V^2_1} +\overline{H^2_2V_2} ) \notag \\
& = {\beta }^2 \big( \ttfrac{1}{8}H_2K^2_3 +\ttfrac{1}{8}H_2(U^2_1+V^2_1) \big) , 
\label{eq-twelvea}
\end{align} 
since the average of $V_1V^2_4$, $K^2_3V_1$, and $H^2_2V_1$ are each $0$. 
\par \noindent Term II is 
\begin{displaymath}
\beta U_4(L_{Y_{V_1}}F) = -{\beta }^2( U_2U_4V_2 + \ttfrac{1}{4} U^2_4V_1). 
\end{displaymath}
So 
\begin{equation}
\beta \overline{U_4(L_{Y_{V_1}}F)} = 0. 
\label{eq-twelveb}
\end{equation}
For term III, we have already shown that 
\begin{equation}
-\ttfrac{\beta }{2}\overline{(L_{Y_F}H_2)K_3} = 0. 
\label{eq-twelvec}
\end{equation}
and for term IV we have already shown that 
\begin{equation}
\ttfrac{\beta }{2}\overline{H_2(L_{Y_{K_3}}F)} = 0. 
\label{eq-twelved}
\end{equation}
Term V is 
\begin{displaymath}
-{\beta }(L_{Y_F}H_2)U_4 = {\beta }^2 \big( U^2_4V_1 + \onehalf H_2K_3U_4 + 
H_2U^2_4 + K_3U_4V_1 \big) . 
\end{displaymath}
So 
\begin{align}
-{\beta }\overline{(L_{Y_F}H_2)U_4} & = {\beta }^2 \big( \overline{U^2_4V_1} + 
\onehalf \overline{H_2K_3U_4} + \overline{H_2U^2_4} + \overline{K_3U_4V_1} \big) \notag \\
& = \ttfrac{{\beta }^2}{2} H_2 (U^2_4 + V^2_4) +\ttfrac{{\beta }^2}{2}K_3(U_4V_1 - U_1V_4), 
\label{eq-twelvee}
\end{align}
since the average of $U^2_4V_1$ and $H_2K_3U_4$ vanish; while 
$\overline{U^2_4} = \onehalf (U^2_4+V^2_4)$ and $\overline{U_4V_1} = U_4V_1 - U_1V_4$. 
\par \noindent Term VI is 
\begin{displaymath}
\beta H_2(L_{Y_{U_4}}F) = {\beta }^2H_2\big( V^2_4 +K^2_3 - \ttfrac{1}{4}K_3U_4 +\ttfrac{1}{4}H_2V_1 
+H^2_2 \big) . 
\end{displaymath}
So 
\begin{align}
\beta \overline{H_2(L_{Y_{U_4}}F)} & = {\beta }^2 \, \overline{V^2_4} +{\beta }^2\, H_2K^2_3 + 
{\beta }^2\, H^3_2 \notag \\
& = \onehalf {\beta }^2\, H_2(U^2_4+V^2_4) +{\beta }^2H_2K^2_3 +{\beta }^2H^3_3,  
\label{eq-twelvef}
\end{align}
since $\overline{K_2U_4} = 0 = \overline{H_2V_1}$. 
\par \noindent Term VII is 
\begin{displaymath}
-\beta (L_{Y_{K_3}}F)V_1 = {\beta }^2\big( \ttfrac{1}{4}( U^2_4V_1-U^2_1V_1-V_1V^2_4 +V^3_1 ) 
+H_2V^2_1 + K_3U_4V_1 \big) . 
\end{displaymath}
So 
\begin{align}
-\beta \overline{(L_{Y_{K_3}}F)V_1} &   
= {\beta }^2 H_2\overline{V^2_1} +{\beta }^2K_3 \overline{U_4V_1} \notag \\
& = \ttfrac{{\beta }^2}{2} H_2(U^2_1+V^2_1) + \ttfrac{{\beta }^2}{2} K_3(U_4V_1 - U_1V_4). 
\label{eq-twelveg}
\end{align}
Term VIII is 
\begin{displaymath}
-\beta K_3(L_{Y_{V_1}}F) = {\beta }^2K_3(U_2V_2 +\ttfrac{1}{4}U_4V_1) . 
\end{displaymath}
So 
\begin{equation}
-\beta \overline{K_3(L_{Y_{V_1}}F)}  = {\beta }^2K_3\overline{U_2V_2} 
+\ttfrac{{\beta }^2}{4} K_3\overline{U_4V_1} 
= \ttfrac{{\beta }^2}{4} K_3(U_4V_1 - U_1V_4), 
\label{eq-twelveh}
\end{equation}
since 
\begin{align*}
\overline{U_2V_2} & = \ttfrac{1}{\pi}\int^{\pi }_0 (U_2\cos 2t +V_2\sin 2t)(-U_2\sin 2t + V_2\cos 2t ) \, \dee t \\
& = \ttfrac{1}{\pi} \int^{\pi }_0 \big( -\onehalf U^2_2 \, \sin 4t +U_2V_2 \, {\cos }^2\, 2t -U_2V_2 \, {\sin }^2\, 2t 
+ \onehalf V^2_2 \, \sin 4t \big) \, \dee t \\
& = \ttfrac{1}{2} (U_2V_2 - U_2V_2) =0. 
\end{align*}
\end{subequations}
Collecting together the results of all the above term calculations gives
\begin{align*}
\beta \overline{L^2_{Y_F}H_2} & = 
\beta \overline{(L_{Y_{U_4}}F)V_1} - \beta \overline{(L_{Y_F}H_2)U_4} 
+{\beta } \overline{H_2(L_{Y_{U_4}}F)} -{\beta }\overline{(L_{Y_{K_3}}F)V_1} 
- \beta \overline{K_3(L_{Y_{V_1}}F)} \\
& \hspace{-.25in} = {\beta }^2\big( [ \ttfrac{1}{8}H_2K^2_3 +\ttfrac{1}{8}H_2(U^2_1+V^2_1) ] 
+ [ \ttfrac{1}{2} H_2(U^2_4+V^2_4) + \ttfrac{1}{2} K_3(U_4V_1 -U_1V_4) ] \\
&  + [ \ttfrac{1}{2}H_2(U^2_4+V^2_4) + H_2K^2_3 + H^3_2 ] 
+ [ \ttfrac{1}{2}H_2(U^2_1+V^2_1) +\ttfrac{1}{2}K_3(U_4V_1-U_1V_4)] \\
& \hspace{.25in} + \ttfrac{1}{4} K_3(U_4V_1 - U_1V_4) \big) \\
&\hspace{-.25in} = {\beta }^2 \big[ (\ttfrac{9}{8}H_2K^2_3 +H^3_2) +\ttfrac{5}{8} H_2(U^2_1+V^2_1) +
H_2(U^2_4+V^2_4) +\ttfrac{5}{4}K_3(U_4V_1-U_1V_4) \big] . 
\end{align*}%
Thus the second order normal form of the regularized Stark Hamiltonian $\mathcal{H}$ on 
${\Xi }^{-1}(0)/S^1$ is 
\begin{align}
{\mathcal{H}}^{(2)}_{\mathrm{nf}} & = H_2 -\ttfrac{3}{2}\varepsilon \beta H_2K_3 - 
\onehalf {\varepsilon }^2 \beta \, \overline{L^2_{Y_F}H_2} \notag \\
& = H_2 -\ttfrac{3}{2}\varepsilon \beta H_2K_3 -\onehalf {\varepsilon }^2 {\beta }^2 
\big[ \ttfrac{9}{8} H_2K^2_3 +H^3_2 \notag \\
& \hspace{.25in} + \ttfrac{5}{8}H_2(U^2_1+V^2_1) +H_2(U^2_4+V^2_4) 
+\ttfrac{5}{4}K_3(U_4V_1-U_1V_4) \big] .
\label{eq-thirteen}
\end{align}

\section{The first order normal form of ${\mathcal{H}}^{(2)}_{\mathrm{nf}}$ on $T_hS^3_1$}

Since $L_{X_{H_2}}{\mathcal{H}}^{(2)}_{\mathrm{nf}} =0$ by construction, the second order 
normal form ${\mathcal{H}}^{(2)}_{\mathrm{nf}}$ (\ref{eq-thirteen}) is a smooth function on 
$(H^{-1}_2(h) \cap {\Xi}^{-1}(0))/S^1= T_hS^3_1$, the tangent $h$-sphere bundle of 
the unit $3$-sphere $S^3_1$, given by  
\begin{align}
{\mathcal{H}}^{(2)}_{\mathrm{nf}} & = h -\ttfrac{3}{2}\varepsilon \beta h K_3 - 
\onehalf {\varepsilon }^2{\beta }^2\big[ \ttfrac{9}{8} K^2_3 + h^3 +\ttfrac{5}{8}(U^2_1+V^2_1) +h(U^2_4+V^2_4) \notag \\
& \hspace{.5in} + K_3(U_4V_1-U_2V_4) \big] \notag \\
& = -\ttfrac{3}{2}\varepsilon \beta h K_3 -\onehalf {\varepsilon }^2{\beta }^2 \big[ \ttfrac{9}{8} K^2_3 +
\ttfrac{5}{8}(U^2_1+V^2_1) +h(U^2_4+V^2_4) + K_3(U_4V_1-U_1V_4) \big] , \notag  
\end{align}
dropping the constant $h-\onehalf {\varepsilon }^{2}{\beta }^2 h^3$. Rescaling time 
$\dee t \mapsto -\frac{3}{2}\varepsilon \beta \dee s$ gives the Hamiltonian 
\begin{equation}
\widetilde{\mathcal{H}} = h K_3 + \ttfrac{1}{3} \varepsilon \beta \big[ \ttfrac{9}{8} K^2_3 +
\ttfrac{5}{8}(U^2_1+V^2_1) + K_3(U_4V_1-U_1V_4) \big] 
\label{eq-fourteen}
\end{equation}
on $T_hS^3_1$.  \medskip 

We now show that the Hamiltonian $\widetilde{\mathcal{H}}$ (\ref{eq-fourteen}) on $T_hS^3_1$ can 
be normalized again. On $(T{\R }^4, {\omega }_4)$ the Hamiltonian 
\begin{displaymath}
K_3(q,p) = \onehalf (q^2_3+q^2_4+p^2_3+p^2_4-q^2_1-q^2_2-p^2_1-p^2_2) 
\end{displaymath}
gives rise to the Hamiltonian vector field $X_{K_3}$, whose flow is 
\begin{align*}
{\varphi }^{X_{K_3}}_t(q,p) & =  \\
&\hspace{-.5in}= \big( q_1 \cos t -p_1\sin t, q_2\cos t -p_2 \sin t, q_3 \cos t + p_3 \sin t , 
q_4 \cos t +p_4 \sin t ,\\
&\hspace{-.25in}q_1 \sin t +p_1 \cos t, q_2 \sin t +p_2 \cos t, -q_3\sin t +p_3 \cos t, 
-q_4 \sin t + p_4\cos t  \big) ,  
\end{align*}
which is periodic of period $2\pi $.  \medskip 

The vector field $X_{K_3}$ on $T{\R }^4$ induces the vector field 
\begin{align*}
Y_{K_3} & = -2L_2 \frac{\partial }{\partial K_1} +2L_1 \frac{\partial }{\partial K_2} 
-2K_2 \frac{\partial }{\partial L_1} +2K_1\frac{\partial }{\partial L_2 } \\
& \hspace{.5in} - 2U_4\frac{\partial }{\partial U_1} + 2U_1 \frac{\partial }{\partial U_4} 
-2V_4\frac{\partial }{\partial V_1} + 2V_1 \frac{\partial }{\partial V_4} , 
\end{align*}
on ${\Xi}^{-1}(0)/S^1 \subseteq {\R }^{16}$ with coordinates $(K,L, H_2, \Xi ; U, V)$, whose flow  
\clearpage
\begin{align*}
{\varphi }^{Y_{K_3}}_s(K,L,H_2, \Xi ; U, V) & =  
\big( -L_2 \sin 2s + K_1 \cos 2s, L_1 \sin 2s +K_2 \cos 2s,  K_3 \\
& \hspace{-1.35in} -K_2\sin 2s +L_1\cos 2s , K_1 \sin 2s + L_2 \cos 2s, L_3, H_2, \Xi ; 
U_1 \cos 2s -U_4 \sin 2s,  \\
&\hspace{-1.25in}U_2, U_3, U_1 \sin 2s +U_4 \cos 2s, 
V_1 \cos 2s - V_4 \sin 2s, V_2, V_3, V_1\sin 2s + V_4 \cos 2s \big) 
\end{align*}
is periodic of period $\pi $. Since $L_{Y_{K_3}}$ maps the ideal of smooth functions which 
vanish identically on ${\Xi}^{-1}(0)/S^1$ into itself, $Y_{K_3}$ is a vector field on ${\Xi}^{-1}(0)/S^1$. 
Since $L_{X_{K_3}}H_2 =0$, it follows that $Y_{K_3}$ induces a vector field on $T_hS^3_1$ with 
periodic flow. So we can normalize again. \medskip 

To compute the normal form of the Hamiltonian 
$\widetilde{\mathcal{H}}$ (\ref{eq-fourteen}) on $T_hS^3_1$ we need only calculate 
the average of the term
\begin{displaymath}
T = \ttfrac{9}{8} K^2_3 + \ttfrac{5}{8}(U^2_1+V^2_1) +h(U^2_4+V^2_4) + K_3(U_4V_1-U_1V_4) 
\end{displaymath}
over the flow ${\varphi }^{Y_{K_3}}_s$. Since $L_{Y_{K_3}}K_3 =0$ and 
\begin{displaymath}
L_{Y_{K_3}}(U_4V_1-U_1V_4) = 2U_1V_1-2U_4V_4+2U_4V_4-2U_1V_1 =0, 
\end{displaymath}
we need only calculate $\overline{U^2_1}$, $\overline{U^2_4}$, $\overline{V^2_1}$, and 
$\overline{V^2_4}$. Now 
\begin{align*}
\overline{U^2_1} & = \ttfrac{1}{\pi} \int^{\pi }_0 (U_1 \cos 2s -U_4 \sin 2s)^2 \, \dee s \\
& = \ttfrac{1}{\pi} \int^{\pi }_0 \big( U^2_1 {\cos }^2\, 2s - U_1U_4\sin 4s + U^2_4 {\sin }^2\, 2s \big) \dee s 
= \onehalf (U^2_1+U^2_4). 
\end{align*}
Similarly, $\overline{V^2_1} = \onehalf (V^2_1+V^2_4)$, $\overline{U^2_4} = \onehalf (U^2_1+U^2_4)$, 
and $\overline{V^2_4} = \onehalf (V^2_1+V^2_4)$. Thus 
\begin{align}
\overline{T} & = (\ttfrac{9}{8}-h)K^2_3 +\ttfrac{13}{16}h(U^2_1+V^2_1 +U^2_4+V^2_4) .
\label{eq-fifteen}
\end{align}
So the first order normal form of $\widetilde{\mathcal{H}}$ (\ref{eq-fourteen}) on $T_hS^3_1$ is 
\begin{equation}
{\widetilde{\mathcal{H}}}^{(1)}_{\mathrm{nf}} = hK_3 +\ttfrac{1}{3}\varepsilon \beta 
[ (\ttfrac{9}{8} - h)K^2_3 + \ttfrac{13}{16}(U^2_1+V^2_1+U^2_2+V^2_2)] 
\label{eq-fifteenstar}
\end{equation}

\section{The reduced Hamiltonian ${\widetilde{\mathcal{H}}}^{(1)}_{\mathrm{nf}}$ on $S^2_h\times S^2_h$}

The polynomial $U^2_1+V^2_1 +U^2_4+V^2_4$ is invariant under 
the flows ${\varphi }^{X_{H_2}}_t$, ${\varphi }^{X_{\Xi }}_u$, and thus is a polynomial on the orbit space $T_hS^3_1/S^1 = S^2_h \times S^2_h$, defined by 
\begin{align*}
K^2_1+K^2_2+K^2_3 + L^2_1+L^2_2+L^2_3 & = h \\
K_1L_1 +K_2L_2 +K_3L_3 & = 0. 
\end{align*}
We now find this polynomial. From the explicit description of ${\R}^8/S^1$ in (\ref{eq-fivenwstar}) it follows that on $\big( H^{-1}_2(h) \cap {\Xi}^{-1}(0) \big) /S^1$ 
\clearpage
\begin{align*}
U_2V_1 - U_1V_2 & = -h K_1 \\
U_3V_1-U_1V_3 & = -hK_2 \\
U_4V_1 -U_2V_4 & =-h K_3.
\end{align*}
So on $S^2_h \times S^2_h$ 
\begin{align*}
h^2(K^2_1+K^2_2+K^2_3) & = (U_2V_1 - U_1V_2)^2 +(U_3V_1-U_1V_3)^2 +(U_4V_1-U_1V_4)^2 \\
& = (U^2_1+U^2_2+U^2_3+U^2_4)V^2_1 - U^2_1V^2_1 \\
& \hspace{.25in} -2(U_1V_1)(U_1V_1+U_2V_2+U_3V_3+U_4V_4) + 2U^2_1V^2_1 \\
& \hspace{.5in} +(V^2_1+V^2_2+V^2_3+V^2_4)U^2_1 - U^2_1V^2_1 \\
& = h^2(V^2_1+U^2_1), 
\end{align*}
since 
$\langle U,U \rangle  =h^2$, 
$\langle V,V \rangle  = h^2$, and 
$\langle U,V \rangle =0$. Thus $U^2_1+V^2_1 = K^2_1+K^2_2+K^2_3$. Again from the explicit description of $\big( H^{-1}_2(h) \cap {\Xi }^{-1}(0)\big) /S^1$ we have 
\begin{align*}
U_4V_3 -U_3V_4 & = -hL_1 \\
U_4V_2 - U_2V_4 & = h L_2 \\
U_4V_1-U_1V_4 &= -hK_3. 
\end{align*}
So on $S^2_h \times S^2_h$ 
\begin{align*}
h^2(L^2_1+L^2_2+K^2_3) & = (U_4V_3 - U_3V_4)^2 +(U_4V_2-U_2V_4)^2 +(U_4V_1-U_1V_4)^2 \\
& = (V^2_1+V^2_2+V^2_3+V^2_4)U^2_4 - U^2_4V^2_4 \\
& \hspace{.25in} -2(U_4V_4)(U_1V_1+U_2V_2+U_3V_3+U_4V_4) + 2U^2_4V^2_4 \\
& \hspace{.5in} +(U^2_1+U^2_2+U^2_3+U ^2_4)V^2_4 - U^2_4V^2_4 \\
& = h^2(U^2_4 +V^2_4). 
\end{align*}
Thus $U^2_4+V^2_4 = L^2_1+L^2_2+K^2_3$. Consequently  
\begin{displaymath}
U^2_1+V^2_1+U^2_4+V^2_4 = K^2_1+K^2_2 +2K^2_3 +L^2_1+L^2_2 
\end{displaymath}
$S^2_h \times S^2_h$. Hence on $S^2_h \times S^2_h$  
\begin{align}
\widehat{\mathcal{H}} = {\widetilde{\mathcal{H}}}^{(1)}_{\mathrm{nf}} & = 
h K_3 +\ttfrac{1}{3}\varepsilon \beta [\ttfrac{1}{16}(18 - 3h) K^2_3 + \ttfrac{13}{16}h (K^2_1 + K^2_2 +K^2_3+L^2_1 +L^2_2)] \notag \\
& = h K_3 +\ttfrac{1}{3}\varepsilon \beta [\ttfrac{1}{16}(18 - 3h) K^2_3 + \ttfrac{13}{16}h (h^2 - L^2_3)].  
\label{eq-sixteen}
\end{align}

\section{The Hamiltonian system $(\widehat{\mathcal{H}}, S^2_h \times S^2_h, \{ \, \, , \, \, \} )$}

Using the coordinates $(\xi , \eta ) = \big( (K+L)/2, (K-L)/2 \big) $ on 
${\R }^3 \times {\R }^3$, the space of smooth functions on the reduced space $S^2_h \times S^2_h$, defined by 
\begin{displaymath}
{\xi }^2_1+{\xi }^2_2+{\xi }^2_3  = h^2\, \, \, \mathrm{and} \, \, \, {\eta }^2_1+{\eta }^2_2 + {\eta }^2_3  = h^2, 
\end{displaymath}
has a Poisson structure with bracket relations
\begin{displaymath}
\{ {\xi }_i, {\xi }_j \} = \sum^3_{k=1}{\epsilon }_{ijk}{\xi }_k, \, \, \, 
\{ {\eta }_i, {\eta }_j \} = -\sum^3_{k=1}{\epsilon }_{ijk}{\eta }_k, \, \, \,
\{ {\xi }_i, {\eta }_j \} =0.
\end{displaymath}
Since $\{ K_3, L_3 \} =0$, it follows that $\{ K_3, \widehat{\mathcal{H}} \} =0$. Thus the flow 
${\varphi }^{Z_{K_3}}_r$ of the Hamiltonian vector field $Z_{K_3}$ on $(S^2_h \times S^2_h, 
\{ \, \, , \, \, \} )$ generates an $S^1$ symmetry of the Hamiltonian system 
$(\widehat{\mathcal{H}}, S^2_h \times S^2_h, \{ \, \, , \, \, \} )$. So this system is 
completely integrable. \medskip 

We reduce this $S^1$ symmetry as follows. Consider the vector field $Z_{K_3}$ on 
${\R }^3 \times {\R }^3$ corresponding 
to the Hamiltonian $K_3 = \onehalf ({\xi }_3 + {\eta }_3)$. Its integral curves satisfy
\begin{align*}
{\dot{\xi}}_1 & = \{ {\xi }_1, K_3 \} = \onehalf \{ {\xi }_1, {\xi }_3 \} = - \onehalf {\xi }_2 \\
{\dot{\xi}}_2 & = \{ {\xi }_2, K_3 \} = \onehalf \{ {\xi }_2, {\xi }_3 \} =  \onehalf {\xi }_1 \\
{\dot{\xi}}_3 & = \{ {\xi }_3, K_3 \} = 0 \\
{\dot{\eta}}_1 & = \{ {\eta }_1, K_3 \} = \onehalf \{ {\eta }_1, {\eta }_3 \} =  \onehalf {\eta }_2 \\
{\dot{\eta}}_2 & = \{ {\eta }_2, K_3 \} = \onehalf \{ {\eta }_2, {\eta }_3 \} =  -\onehalf {\eta }_1 \\
{\dot{\xi}}_3 & = \{ {\eta }_3, K_3 \} = 0.
\end{align*}
Thus the flow of $Z_{K_3}$ on ${\R }^3 \times {\R }^3$ is 
\begin{align*}
{\varphi }^{Z_{K_3}}_t(\xi , \eta ) & = 
\big( {\xi }_1 \cos t/2 -{\xi}_2 \sin t/2, {\xi }_1 \sin t/2 +{\xi }_2 \cos t/2, {\xi }_3 , \\
& \hspace{.5in} {\eta }_1 \cos t/2 +{\eta }_2 \sin t/2, {\eta }_1 \sin t/2 -{\eta }_2 \cos t/2, {\eta }_3 \big) , 
\end{align*}
which preserves $S^2_h \times S^2_h$ and is periodic of period $4\pi $. \medskip 

We now determine the space $(S^2_h\times S^2_h)/S^1$ of orbits of the vector field 
$Z_{K_3}$. We use invariant theory. The algebra of polynomials on ${\R }^3 \times {\R }^3$, 
which are invariant under the $S^1$ action given by the flow ${\varphi }^{Z_{K_3}}_t$, is 
generated by 
\begin{displaymath}
\begin{array}{lclcl}
{\sigma }_1 = {\xi }^2_1 +{\xi }^2_2 & & {\sigma }_2 = {\eta }^2_1 +{\eta }^2_2 & & 
{\sigma }_3 = {\xi }_1{\eta }_2 - {\xi}_2{\eta }_1 \\
\rule{0pt}{12pt} {\sigma }_4 = {\xi }_1{\eta }_1 +{\xi }_2{\eta }_2 & & {\sigma }_5 = \onehalf ({\xi }_3+{\eta }_3) & & {\sigma }_6 = \onehalf ({\xi }_3-{\eta }_3) ,  
\end{array} 
\end{displaymath} 
which are subject to the relation 
\begin{subequations}
\begin{align}
{\sigma }^2_3 +{\sigma }^2_4 & = ({\xi }_1{\eta }_2-{\xi}_2{\eta }_1)^2 + ({\xi }_1{\eta }_1+{\xi}_2{\eta }_2)^2 
\notag \\
& = ({\xi }^2_1+{\xi }^2_2)({\eta }^2_1+{\eta }^2_2) = {\sigma }_1 {\sigma }_2, \, \, \, 
{\sigma }_1 \ge 0 \, \, \& \, \, {\sigma }_2 \ge 0. 
\label{eq-nineteena}
\end{align}
In terms of invariants the defining equations of $S^2_h \times S^2_h$ become
\begin{align}
{\sigma }_1 +({\sigma }_5+{\sigma }_6)^2 & = {\xi }^2_1+{\xi }^2_2 +{\xi }^2_3 = h^2 
\label{eq-nineteenb} \\
{\sigma }_2 +({\sigma }_5-{\sigma }_6)^2 & = {\eta }^2_1+{\eta }^2_2 +{\eta }^2_3 = h^2 
\label{eq-nineteenc}
\end{align}
\end{subequations}
Eliminating ${\sigma }_1$ and ${\sigma }_2$ from (\ref{eq-nineteena}) using (\ref{eq-nineteenb}) 
and (\ref{eq-nineteenc}) gives 
\begin{subequations}
\begin{equation}
{\sigma }^2_3 +{\sigma }^2_4 = \big( h^2-({\sigma }_5+{\sigma }_6)^2 \big) 
\big( h^2 - ({\sigma }_5 - {\sigma }_6)^2 \big) , \, \, |{\sigma }_5 +{\sigma }_6| \le h \, \, \& \, \, 
|{\sigma }_5-{\sigma }_6| \le h, 
\label{eq-twentya}
\end{equation}
which defines $(S^2_h \times S^2_h)/S^1$ as a semialgebraic variety in ${\R }^4$ with 
coordinates $({\sigma}_3, {\sigma }_4, {\sigma }_5, {\sigma }_6)$. Thus the reduced space 
$\big( K^{-1}_3(2 k ) \cap (S^2_h\times S^2_h) \big)/S^1$ is defined by (\ref{eq-twentya}) and 
\begin{equation}
{\sigma }_5 = k . 
\label{eq-twentyb}
\end{equation}
\end{subequations}
Consequently, $\big( K^{-1}_3 \cap (S^2_h \times S^2_h) \big) (2 \ell )/S^1$ is the semialgebraic variety 
\begin{align}
{\sigma }^2_3+{\sigma }^2_4 & = \big( h^2-(k +{\sigma }_6)^2 \big) 
\big( h^2-(k-{\sigma }_6)^2 \big) \notag \\
& = \big( (h-k )^2-{\sigma }^2_6 \big) \big( (h+k )^2-{\sigma }^2_6\big) , \, \, 
|{\sigma }_6| \le h - |\ell |  
\label{eq-twentyone}
\end{align} 
in ${\R }^3$ with coordinates $({\sigma }_3, {\sigma }_4, {\sigma }_6)$. When 
$0 < |k | < h$ the reduced space (\ref{eq-twentyone}) is a smooth $2$-sphere. 
When $|k| = h$ it is a point. When $k =0$ it is a topological $2$-sphere with 
conical singular points at $(0,0,\pm h)$. These singular points correspond to the 
fixed points $h(0,0, \pm 1, 0, 0, \mp 1)$ of the $S^1$ action on $S^2_h \times S^2_h$ 
generated by the flow of the vector field $Z_{K_3}$. \medskip 

By (\ref{eq-sixteen}) the reduced Hamiltonian on 
$\big( K^{-1}_3(2k) \cap (S^2_h \times S^2_h) \big)/S^1$ is 
\begin{equation}
{\widehat{\mathcal{H}}}_{\mathrm{red}} = -\ttfrac{13}{12}\varepsilon \beta {\sigma }^2_6, 
\label{eq-twentytwo}
\end{equation}
using $L_3 = {\xi }_3-{\eta }_3 = 2{\sigma }_6$, having dropped the constant 
$hk + \ttfrac{1}{48} \varepsilon \beta [ (18 -3h)k^2 + 13 h^3] $.

\end{document}